\documentclass[11pt]{article}

\usepackage{amssymb}

\def\A{{\cal A}}
\def\P{{\cal P}}
\def\C{{\cal C}_{\mathbb Q}}
\def\N{{\mathbb N}}
\def\Z{{\mathbb Z}}
\def\RR{{\mathbb R}}
\def\R{{\cal R}}
\def\QQ{{\mathbb Q}}
\def\Q{{\mathbb Q}^*_+}

\newenvironment{verif}{\noindent{\it Proof.}}
{\newline\mbox{\ }\hfill\rule{0.5em}{0.5em}\smallskip}

\begin{document}

\title{\bf Ergodicity of the action of the positive
rationals on the group of  finite adeles and the Bost-Connes phase
transition theorem}

\author{Sergey Neshveyev}

\date{}

\maketitle

\begin{abstract}
For each $\beta\in(0,+\infty)$ there exists a canonical measure
$\mu_\beta$  on the ring $\A_f$ of finite adeles. We show that
$\Q$ acts ergodically on $(\A_f,\mu_\beta)$ for $\beta\in(0,1]$,
and then deduce from this the uniqueness of KMS$_\beta$-states for
the Bost-Connes system.
\end{abstract}

\bigskip\bigskip

Bost and Connes \cite{BC} constructed a remarkable C$^*$-dynamical
system which has a phase transition with spontaneous symmetry
breaking involving an action of the Galois group
Gal($\QQ^{ab}/\QQ$), and whose partition function is the Riemann
$\zeta$ function. In their original definition the underlying
algebra arises as the Hecke algebra associated with an inclusion
of certain $ax+b$ groups. Recently Laca and Raeburn~\cite{LR,L2}
have realized the Bost-Connes algebra as a full corner of the
crossed product algebra $C_0(\A_f)\rtimes\Q$. This new look at the
system has allowed to simplify significantly the proof of the
existence of KMS-states for all temperatures, and the proof of the
phase transition theorem for $\beta>1$~\cite{L1}. On the other
hand, for $\beta\le1$ the uniqueness of KMS$_\beta$-states implies
the ergodicity of the action of $\Q$ on $\A_f$ for certain
measures (in particular, for the Haar measure). The aim of this
note is to give a direct proof of the ergodicity, and then to show
that the uniqueness of KMS$_\beta$-states easily follows from it.
Though the proof of the Bost-Connes phase transition theorem (for
$\beta\le1$) thus obtained differs from the proofs given
in~\cite{BC} and~\cite{L1}, it is entirely based on these papers.
In particular, the key point is an application of Dirichlet's
theorem.

\bigskip

So let $\P$ be the set of prime numbers, $\A_f$ the restricted
product of the fields $\QQ_p$, $p\in\P$, of $p$-adic numbers,
$\R=\prod_p\Z_p$ its maximal compact subring,
$W=\R^*=\prod_p\Z^*_p$. The group $\Q$ of positive rationals is
embedded diagonally into $\A_f$, and so acts by multiplications on
the additive group of finite adeles. Then the Bost-Connes algebra
$\C$ is the full corner of $C_0(\A_f)\rtimes\Q$ determined by the
characteristic function of $\R$~\cite{L2}. The dynamics $\sigma_t$
is defined as follows~\cite{L1}: it is trivial on $C_0(\A_f)$, and
$\sigma_t(\lambda(q))=q^{it}\lambda(q)$, where $\lambda(q)$ is the
multiplier of $C_0(\A_f)\rtimes\Q$ corresponding to $q\in\Q$. Then
(\cite{L1}) there is a one-to-one correspondence  between
$(\beta,\sigma_t)$-KMS-states on $\C$ and measures $\mu$ on $\A_f$
such that
$$
\mu(\R)=1\ \ \hbox{and}\ \ q_*\mu=q^\beta\mu\ \ \hbox{for all}\
q\in\Q\ \ (\hbox{i.e.,}\ \mu(q^{-1}X)=q^\beta\mu(X)).\eqno(1\beta)
$$
Namely, the KMS-state corresponding to $\mu$ is the restriction of
the dual weight on $C_0(\A_f)\rtimes\Q$ to $\C$.

Note that if $\beta>1$ and $\mu$ is a measure with the property
(1$\beta$) then $\mu(W)={1\over\zeta(\beta)}>0$, since
$W=\R\backslash\cup_pp\R$. Moreover, the sets $qW$, $q\in\Q$, are
disjoint, and their union is a set of full measure (since
$\sum_{n\in\N}\mu(nW)=1$). Thus there exists a one-to-one
correspondence between probability measures on $W$ and measures on
$\A_f$ satisfying (1$\beta$)~\cite{L1}. On the other hand, if
$\beta\le1$ then $\mu(W)=0$.

For each $\beta\in(0,+\infty)$ there is a unique $W$-invariant
measure $\mu_\beta$ satisfying (1$\beta$)~\cite{BC,L1}.
Explicitly, $\mu_\beta=\otimes_p\mu_{\beta,p}$, where
$\mu_{\beta,p}$ is the measure on $\QQ_p$ such that $\mu_{1,p}$ is
the Haar measure ($\mu_{1,p}(\Z_p)=1$), and
$$
{d\mu_{\beta,p}\over
d\mu_{1,p}}(a)={1-p^{-\beta}\over1-p^{-1}}|a|^{\beta-1}_p\ \
\hbox{for}\ \ a\in\QQ_p.
$$
In fact, for the proof of Proposition below we will only need to
know that the restriction of $\mu_{\beta,p}$ to $\Z^*_p$ is a
(non-normalized) Haar measure.

\medskip\noindent
{\sc Proposition.} {\it The action of $\Q$ on $(\A_f,\mu_\beta)$
is ergodic for $\beta\in(0,1]$.}

\medskip
\begin{verif}
Consider the space $L^2(\R,d\mu_\beta)$ and the subspace $H$ of it
consisting of the functions that are constant on $\N$-orbits. In
other words, $H=\{f\in L^2(\R,d\mu_\beta)\,|\,V_nf=f,\ n\in\N\}$,
where $(V_nf)(x)=f(nx)$. Since any $\Q$-invariant subset of $\A_f$
is completely determined by its intersection with $\R$, it
suffices to prove that $H$ consists of constant functions. For
this we will compute the action of the projection $P\colon
L^2(\R,d\mu_\beta)\to H$ on a basis of $L^2(\prod_{p\in
B}\Z_p,\otimes_{p\in B}\mu_{\beta,p})$ (considered as a subspace
of $L^2(\R,d\mu_\beta)$) for each finite subset $B$ of $\P$.

Let $\chi$ be a character of $\prod_{p\in B}\Z^*_p$. Consider
$\chi$ first as a function on $\prod_{p\in B}\Z_p$ by letting
$\chi=0$ outside of $\prod_{p\in B}\Z^*_p$. Then using the
projection $\R\to\prod_{p\in B}\Z_p$, consider $\chi$ as a
function on $\R$. Let $\N_B$ be the unital multiplicative
subsemigroup of $\N$ generated by $p\in B$. Note that the sets
$n\prod_{p\in B}\Z^*_p$, $n\in\N_B$, are disjoint, their union is
a subset of $\prod_{p\in B}\Z_p$ of full measure, and the operator
$n^{-\beta/2}V^*_n$ maps isometrically $L^2(\prod_{p\in
B}\Z^*_p,\otimes_{p\in B}\mu_{\beta,p})$ onto $L^2(n\prod_{p\in
B}\Z^*_p,\otimes_{p\in B}\mu_{\beta,p})$ for any $n\in\N_B$. Hence
the functions $V^*_n\chi$, $n\in\N_B$, $\chi\in(\prod_{p\in
B}\Z^*_p)\hat{}$, form an orthogonal basis for $L^2(\prod_{p\in
B}\Z_p,\otimes_{p\in B}\mu_{\beta,p})$. So we have to compute
$PV^*_n\chi$. But if $g\in H$ then $(V^*_n\chi,g)=(\chi,g)$,
whence $PV^*_n\chi=P\chi$. Thus we have only to compute $P\chi$.

For a finite subset $A$ of $\P$, let $H_A$ be the subspace
consisting of the functions that are constant on $\N_A$-orbits,
$P_A$ the projection onto $H_A$. Then $P_A\searrow P$ as
$A\nearrow\P$. Set
$$
W_A=\prod_{p\in A}\Z^*_p\times\prod_{q\in\P\backslash
A}\Z_q\subset\R.
$$
Note, as above, that $\cup_{n\in\N_A}nW_A$ is a subset of $\R$ of
full measure. We state that
$$
P_Af|_{\N_Ax}\equiv{1\over\zeta_A(\beta)}\sum_{n\in\N_A}n^{-\beta}f(nx)\
\ \hbox{for}\ \ x\in W_A, \eqno(2)
$$
where $\zeta_A(\beta)=\sum_{n\in\N_A}n^{-\beta}=\prod_{p\in
A}(1-p^{-\beta})^{-1}$. Indeed, denoting the right hand part of
(2) by~$f_A$, for $g\in H_A$ we obtain
\begin{eqnarray*}
(f_A,g)&=&\sum_{n\in\N_A}\int_{nW_A}f_A(x)\overline{g(x)}d\mu_\beta(x)
    =\sum_{n\in\N_A}n^{-\beta}\int_{W_A}f_A(x)\overline{g(x)}d\mu_\beta(x)\\
&=&\zeta_A(\beta)\int_{W_A}f_A(x)\overline{g(x)}d\mu_\beta(x)
    =\sum_{n\in\N_A}n^{-\beta}\int_{W_A}f(nx)\overline{g(x)}d\mu_\beta(x)\\
&=&\sum_{n\in\N_A}\int_{nW_A}f(x)\overline{g(x)}d\mu_\beta(x)=(f,g).
\end{eqnarray*}

Returning to the computation of $P\chi$, we see that
$$
P_A\chi|_{\N_Ax}\equiv{\chi(x)\over\zeta_A(\beta)}
\sum_{n\in\N_A}n^{-\beta}\chi(n) =\chi(x)\prod_{p\in
A}{1-p^{-\beta}\over1-\chi(p)p^{-\beta}}\ \ \hbox{for}\ \ x\in
W_A.
$$
Thus if $\chi$ is trivial then $P_A\chi\equiv\prod_{p\in
B}(1-p^{-\beta})$ for all $A\supset B$, hence $P\chi$ is a
constant. If $\chi$ is non-trivial then since
$||P_A\chi||_\infty\le1$ and the product
$\prod_{p:\,Re\chi(p)<0}(1-p^{-\beta})$ diverges by Dirichlet's
theorem~\cite{S}, we have $P\chi=0$. \end{verif}

\medskip\noindent
{\sc Corollary.}\cite{BC} {\it For $\beta\in(0,1]$ there exists a
unique $(\beta,\sigma_t)$-KMS state on $\C$.}

\medskip
\begin{verif}
Let $\phi_\beta$ be the KMS$_\beta$-state corresponding to
$\mu_\beta$. Since $L^\infty(\A_f,d\mu_\beta)\rtimes\Q$ is a
factor by Proposition, and $\pi_{\phi_\beta}(\C)''$ is its
reduction, $\phi_\beta$ is a factor state. This and the discussion
before Proposition show that
\begin{list}{}{}
\item{(i)} $\phi_\beta$ is an extremal KMS$_\beta$-state;
\item{(ii)} $\phi_\beta$ is a unique $W$-invariant KMS$_\beta$-state.
\end{list}
Now the proof is finished as in~\cite[Theorem 25]{BC}:

If $\psi$ is an extremal KMS$_\beta$-state then
$\int_Ww_*\psi\,dw=\phi_\beta$. Since KMS$_\beta$-states form a
simplex, we conclude that $\psi=\phi_\beta$.
\end{verif}

\medskip\noindent
{\sc Remarks.}

\noindent (i) The expression for $P\chi$ in the proof of
Proposition shows that the divergence of the product
$$
\prod_{p\in\P}\left|{1-p^{-\beta}\over1-\chi(p)p^{-\beta}}\right|
$$
for non-trivial $\chi$ is a necessary condition for the ergodicity
(otherwise $P\chi$ would be a non-zero function, which can not be
constant since $\int_\R P\chi\,d\mu_\beta=\int_\R\chi\,
d\mu_\beta=0$), hence for the uniqueness of KMS$_\beta$-states. So
the appearance of (some form of) Dirichlet's theorem in the proofs
is not an accident.

\smallskip\noindent (ii) By \cite[Theorem 5]{BC}
$\pi_{\phi_\beta}(\C)''$ is a factor of type III$_1$ for
$\beta\in(0,1]$. Then the factor
$L^\infty(\A_f,d\mu_\beta)\rtimes\Q$ is also of type III$_1$.
Hence its smooth flow of weights is trivial, that means that the
action of $\Q$ on $(\RR_+\times\A_f,dt\otimes d\mu_\beta)$ is
ergodic~\cite{CT}. In particular, the spectral subspaces of
$L^\infty(\A_f,d\mu_\beta)$ corresponding to the characters
$q\mapsto q^{it}$ of $\Q$ have to be trivial for all $t\ne0$. But
the projection $P_t$ onto the subspace $\{f\,|\,V_nf=n^{it}f\}$ of
$L^2(\R,d\mu_\beta)$ is computed with the same ease as in the
proof of Proposition:
$$
P_t=s-\lim_{A\nearrow\P}P_{t,A},\ \
(P_{t,A}f)(mx)={m^{it}\over\zeta_A(\beta)}\sum_{n\in\N_A}n^{-\beta-it}f(nx)\
\ \hbox{for}\ \ x\in W_A,\ m\in\N_A.
$$
Thus the product
$$
\prod_{p\in\P}\left|{1-p^{-\beta}\over1-\chi(p)p^{-\beta-it}}\right|
$$
has to be divergent for all $t\ne0$ and all number characters
$\chi$ modulo $m$.

\bigskip

\begin{flushleft}
Institute for Low Temperature Physics \& Engineering\\ Lenin Ave
47\\ Kharkov 310164, Ukraine\\ neshveyev@ilt.kharkov.ua\\
\end{flushleft}

\end{document}